\documentclass[graybox]{svmult}

\usepackage{type1cm}        
\usepackage{makeidx}         
\usepackage{graphicx}        
\usepackage{multicol}        
\usepackage[bottom]{footmisc}
\usepackage{xcolor}

\usepackage{newtxtext}       %
\usepackage[varvw]{newtxmath}       
\usepackage{siunitx}
\usepackage{tikz}
\usepackage[colorlinks,citecolor=blue]{hyperref}
\IfFileExists{.git/gitHeadInfo.gin}{
    \usepackage[pcount,grumpy,mark,markifdirty]{gitinfo2}
}{%
    \usepackage[local,pcount,grumpy,mark,markifdirty]{gitinfo2}
}

\usepackage{xspace}
\newcommand{\feelpp}{Feel\nolinebreak\hspace{-.05em}\raisebox{.4ex}{\tiny\bf +}\nolinebreak\hspace{-.10em}\raisebox{.4ex}{\tiny\bf +}\xspace}

\newcommand{\Fluid}{\mathcal{F}} 
\newcommand{\Alemap}{\mathcal{A}} 
\newcommand{\ALE}{ALE} 
\newcommand{\Vel}{u} 
\newcommand{\Pres}{p} 
\newcommand{\Density}{\rho} 
\newcommand{\tvel}{U} 
\newcommand{\angvel}{\omega} 
\newcommand{\Rmat}{R} 
\newcommand{\Angle}{\theta}
\newcommand{\Inertia}{I} 
\newcommand{\mass}{m} 
\newcommand{\CenterMassi}{x_{CM}^i}
\newcommand{\Solid}{\mathcal{S}} 
\newcommand{\normal}{n} 
\newcommand{\CompDomain}{\Fluid}
\newcommand{\dx}{\, \mathrm{d}x}

\newcommand{\R}{\mathbb{R}}
\newlength{\dhatheight}

\makeindex             

\begin{document}

\title*{Towards a computational framework using finite element methods with Arbitrary Lagrangian-Eulerian approach for swimmers with contact}
\titlerunning{A computational framework for swimmers}
\author{Luca Berti, Laetitia Giraldi, Christophe Prud'homme, Céline Van Landeghem}
\institute{Luca Berti \at CNRS - Institut de Recherche Mathématique Avancée, UMR 7501 Université de Strasbourg  \email{berti@math.unistra.fr}
\and Laetitia Giraldi \at Universit\'e C\^ote d'Azur, INRIA, CNRS, Calisto team, \email{laetitia.giraldi@inria.fr}
\and Christophe Prud'homme \at Institut de Recherche Mathématique Avancée, UMR 7501 Université de Strasbourg - CNRS \email{christophe.prudhomme@math.unistra.fr}
\and Céline Van Landeghem \at Institut de Recherche Mathématique Avancée, UMR 7501 Université de Strasbourg - CNRS \email{c.vanlandeghem@unistra.fr}}
%
%
\maketitle

\abstract{Swimming involves a body's capability to navigate through a fluid by undergoing self-deformations. Typically, fluid dynamics are described by the Navier-Stokes equations, and when integrated with a swimming body, it results in a highly intricate model. This paper introduces a computational framework for simulating the movement of multiple swimmers with various geometries immersed in a Navier-Stokes fluid. The approach relies on the finite element method with an Arbitrary Lagrangian-Eulerian (ALE) framework to handle swimmer displacements. Numerous numerical experiments demonstrate the adaptability of the computational framework across various scenarios.
All the implementations are made using the \feelpp\ finite element library \cite{prudhomme_feelppfeelpp_2023}.}

\section{Introduction}


Swimming involves the capacity of a body to navigate through a fluid by undergoing self-deformations. Typically, the fluid dynamics are described by the Navier-Stokes equations, and when coupled with a swimming body, it yields a highly intricate model. One of the main difficulties lies in the fact that the swimmers deform and move, requiring constant adjustment of the fluid domain and fluid-swimmer interface. Other difficulties come into play when accounting for obstacles, the possible presence of other swimmers, and the modeling of their interactions. These interactions are even more complex because disturbances caused by fluid agitation at one point in the domain can propagate throughout the entire domain.




The problem of simulating the motion of deformable swimmers in a fluid can be addressed through various approaches, depending on the flow regime. In the low Reynolds number limit, when the viscous effects prevail over inertial effects, as it is the case for microscopic bodies in very viscous fluids \cite{huang_life_2013}, it is legitimate to consider that the fluid is governed by the Stokes equations \cite{happel_low_2012,lauga_hydrodynamics_2009}. 



In this case, the hydrodynamic drag can be simplified asymptotically by employing the Resistive Force Theory \cite{gray_propulsion_1955}, in which the hydrodynamic friction is related to the fluid velocity via an 
anisotropic relation 
depending on some parameters that account for the shape of the moving object. This approach does not require the numerical solution of the Stokes equations, but it has a
very limited physical validity, often restricted to the case of very simple boundary-free domains.

Slender body theory improves the Resistive Force Theory
to investigate the dynamics of slender particles in highly viscous flows \cite{cox_motion_1970, batchelor_slender-body_1970, keller_slender-body_1976,johnson_improved_1980}. 
This approach involves solving an integral equation which describes the effects of the fluid on the slender body. Assuming small deformations, the equilibrium of forces and torques can be expressed as a partial differential equation \cite{gadelha_counterbend_2013,  moreau_asymptotic_2018, el_khiyati_steering_2023}, whose solution defines the curve that describes the fiber. However, this method neglects the effects of fiber deformation on the fluid, leading to an approximation in accounting for the effects of the displacement of multiple swimmers.

A numerical method that is widely used in the micro-swimming community is the boundary element method, that relies on the integral formulation of the Stokes equations to determine the fluid velocity field \cite{pozrikidis_boundary_1992}. 
This formulation reduces the computational cost of the dynamics of the swimmer by working on the swimmer's boundary alone, which ensures significant memory savings since the fluid domain is not discretized. 
The boundary integral formulation expresses the fluid velocity at the boundary of the swimmer as a function of fluid stresses and velocities, via the fundamental solutions of the Stokes equations. Since these functions exhibit singular behaviour when the source and evaluation points are close together, a regularisation procedure needs to be employed for the numerical solution. Several solutions are possible: a regularisation of the numerical approximation of the singular integrals using a semi-analytic procedure \cite{huang_notes_1993} or the use of regularised kernels \cite{olson_coupling_2011}.


When inertial effects become more important and the low Reynolds number assumption is no longer valid, 
it is necessary to resort to the solution of the system coupling Navier-Stokes equations with Newton's equations for rigid bodies. 
A commonly used method involves approximating the values of fluid velocity and pressure on a fixed discretization of the domain and representing the swimmers using an implicit function that parameterizes their boundaries. 
These approaches are commonly referred to as Immersed Boundary Methods with level set functions describing the swimmer. 
The challenge of these methods lies in following the swimmer. 
For instance, one may use methods called CutFEM (Cut Finite Element Method), which project the swimmer's boundary onto the fixed mesh \cite{monasse_conservative_2012,bergmann_accurate_2014,bergmann_bioinspired_2016,hansbo_cut_2016,burman_cutfem_2015}.



In this chapter, we discuss an alternative method to compute the displacement of multiple deformable swimmers immersed in a Navier-Stokes fluid, which can be confined within a geometrically complex domain. This method is based on the finite element solution of the swimming problem, using the Arbitrary-Lagrangian-Eulerian description for the motion of the fluid domain.
This approach is popular in fluid-structure interaction problems \cite{chabannes_high-order_2013,pena_high_2012}, and it is based on a formulation of the fluid equations in an intermediate frame between the Eulerian and Lagrangian description.
A characteristic of this method is that the solutions are approximated on a mesh that deforms over time, following the deformation and movements of the swimmers.
As the immersed boundary method, it accounts for the full, non-approximated description of the fluid dynamics.
This numerical method is very versatile and it has been used in a wide variety of contexts and physical problems, and it can be easily extended to cases in which the Newtonian fluid is substituted by a complex fluid, which is often the case of biological media.


The structure of the chapter is the following: in Section \ref{Sec:MathModeling} the fluid, swimmer and collision models are presented. In Section \ref{Sec:NumDiscr}, the numerical discretization of the swimming problem is addressed. Finally, in Section \ref{Sec:Applications}, several numerical applications illustrate this framework.









\section{Mathematical modeling}
\label{Sec:MathModeling}
This section details the mathematical modeling of the coupled fluid-swimmer interaction problem. A model of collision forces is also discussed, allowing to simulate swimmer-swimmer or swimmer-boundary interactions.


\subsection{The fluid model}

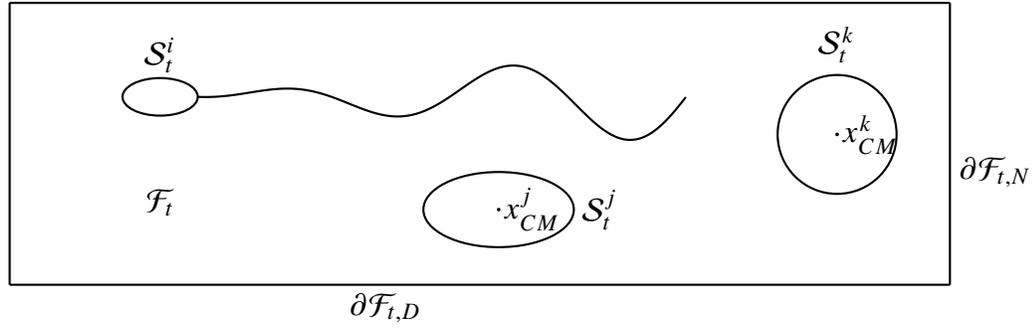
\begin{figure}
    \centering
\begin{center}
\begin{tikzpicture}

\draw[color=black,  thick] (2*5,0*3.75) -- (5*4.5,0)     ;
\draw[color=black,  thick] (2*5,1.*3.75) -- (5*4.5,1.*3.75) ;
\draw[color=black,  thick] (2*5,0) -- (2*5,1.*3.75) ;
\draw[color=black,  thick] (5*4.5,0) -- (5*4.5,1.*3.75) ;

\draw[draw=black,thick] (12,0.5*5) ellipse (0.5cm and 0.25cm);
	\begin{scope}[shift={(-3.0,0.5*5)}] \draw[draw=black,thick] plot[domain=15.5:7*pi, samples=320] (\x,{(\x-15.5)/10*sin(2*\x r)});
    \end{scope} 

\node  at (12,0.61*5) {\large{$\Solid_t^i$}};
        
\draw[color=black, thick] (15.5+1.,1.) ellipse (1*1. and 0.5*1.) ;
\node  at (15.+2.85,1.) {\large{$\Solid_t^j$}};
\draw[color=black, thick] (15.5+1.,1.) circle (0.25*1.1pt);
\node  at (15.5+1.45,1.) {\large{$x_{CM}^j$}};

\draw[color=black, thick](21,2) circle (9*2.5pt);
\node  at (21,3.2) {\large{$\Solid_t^k$}};
\draw[color=black, thick](21,2) circle (0.25*1.1pt);
\node  at (21.45,2) {\large{$x_{CM}^k$}};

\node  at (11.+1.,1.1) {\large{$\Fluid_t$}};

\node  at (23.1,1.5) {\large{$\partial \Fluid_{t,N}$}};
\node  at (15.,-0.3) {\large{$\partial \Fluid_{t,D}$}};

\end{tikzpicture}
\end{center}
   \caption{Notations for the fluid-swimmer model. $\CompDomain_t$ denotes the fluid domain and $\Solid^i_t,\Solid^j_t,\Solid^k_t$ denote the swimmers and solid objects included in the fluid. $\partial \Fluid_{t,N}$ and $\partial \Fluid_{t,D}$ denote the portions of the fluid boundary where Dirichlet and Neumann boundary conditions are prescribed.}
    \label{fig:swimmer-tikz}
\end{figure}

We consider the Navier-Stokes equations in a moving domain to describe the motion of a Newtonian fluid medium. 
Let $\CompDomain_t\subset \R^d$, where $d=2,3$ is the dimension, denote the region occupied by the fluid at time $t$, $\mu$  the constant fluid viscosity and $\Density_\Fluid$ the constant fluid density. Let $\Solid_t^i \subset \R^d$ be the domain that is occupied by the body $i$, $i = 1,...,N$, where $N$ is the number of bodies, of constant density $\Density_\Solid^i$. At each time instant $t$, the fluid is in contact with the swimmers, and $\bar{\CompDomain}_{t} \cap \bar{\Solid}_t^i = \partial \Solid_t^i$.
The notations are summarized in Figure \ref{fig:swimmer-tikz}.

The velocity field $\Vel(t,x):]0,T]\times \CompDomain_t \rightarrow \R^d$ and the pressure field $\Pres(t,x):]0,T]\times\CompDomain_t \rightarrow \R$ satisfy the Navier-Stokes equations given by: 
\begin{equation}
	\begin{aligned}
		\rho_\Fluid \Big(\partial_t \Vel +( \Vel \cdot \nabla )\Vel \Big) - \nabla \cdot \sigma(u,p) &=  f
		\quad &\text{on $ ]0,T]\times \CompDomain_t$},
		\\
		\nabla \cdot \Vel &= 0  \quad&\text{on $ ]0,T]\times \CompDomain_t$},\\
		\Vel &= \bar{\Vel}^i
		\quad &\text{on $]0,T]\times \partial \Solid_t^i$},\\
		\Vel(0,x) &= \Vel_0(x) \quad &\text{on $ \{0\}\times\CompDomain_0$},\\
		\Vel(t,x) &= h(t,x) &\text{on $ ]0,T]\times\partial\CompDomain_{t,D}$},\\
		\sigma(u,p) \normal&=  g(t,x) &\text{on $ ]0,T]\times\partial\CompDomain_{t,N}$},
	\end{aligned}
	\label{Eq:NS}
\end{equation}
where $\sigma(u,p)$ is the stress tensor $\sigma(u,p) = -\Pres\mathbb{I} + 2 \mu D(u)$, with $D(u) = \frac{1}{2}\left(\nabla u + (\nabla u)^T\right)$, and
$\bar{\Vel}^i(t,x):[0,T]\times\partial \Solid_t^i \rightarrow \R^d$ results from the interaction between the fluid and the $i$-th swimmer's body and $f(t,x):[0,T]\times\CompDomain_t \rightarrow \R^d$ represents external volume forces. The function $u_0$ is the prescribed initial condition of the system. 
Here $h(t,x)$ and $g(t,x)$ represent the Dirichlet and Neumann boundary conditions, prescribed on the respective portions of the fluid boundary $ ]0,T]\times\partial\CompDomain_{t,D}$ and  $ ]0,T]\times\partial\CompDomain_{t,N}$, and we have that $\partial \CompDomain_{t} = \partial\CompDomain_{t,D} \cup \partial\CompDomain_{t,N} \cup \partial \Solid_t^i$, $i=1,...,N$, for all times. 

In the case of swimming, $\bar{\Vel}^i$ can be split into two contributions. The first one is determined by the linear velocity $\tvel^i(t): [0,T]\to\R^d$ and angular velocity  $\angvel^i(t):[0,T]\to \R^{d^*}$ of the swimmer around its center of mass $\CenterMassi$, where $d^*=1$ if $d=2$ and $d^*=3$ if $d=3$.  The second one is given by the deformation velocity $u_d^i$, defined on the boundary of the swimmer. This latter models the swimming gait for active particles and vanishes for rigid bodies. The function $\bar{u}^i$ then becomes $\bar{\Vel}^i(t,x) = \tvel^i(t) + \angvel^i(t)\times (x-\CenterMassi)+ \Vel_d^i(t)$.

\subsection{The swimmer model}


The swimmer's motion is either a result of external forces such as gravity and collision forces, or of its deformation, which generates hydrodynamic forces that are then translated into rigid movement through Newton's laws.

Let us define the mass of the $i$-th swimmer by  $\mass^i = \int_{\Solid^i} \Density_{\Solid}^i$, where $\Solid^i = \Solid^i_0$ is the reference configuration, 
and $\Inertia^i = \int_{\Solid^i} \Density_{\Solid}^i (x-\CenterMassi) \otimes (x-\CenterMassi)$ 
its positive definite and symmetric inertia tensor. 
The Newton and Euler equations describing the rigid body velocities $\tvel^i$ and $\angvel^i$ read:
\begin{equation}
	\begin{aligned}
		\mass^i \frac{d}{dt}\tvel^i &= F_e^i -\int_{\partial \Solid^i} \sigma \normal \,\textrm{ds},\\
		\frac{d}{dt}(\Rmat \Inertia^i \Rmat^T \angvel^i) &= T_e^i -\int_{\partial \Solid^i} \sigma \times (x-\CenterMassi)\,\textrm{ds}, 
	\end{aligned}
	\label{Eq:RB}
\end{equation}

where $\normal$ is the unit outward normal to $\partial \Solid^i$. Equations \eqref{Eq:RB} correspond to the balance of forces and torques applied to each swimmer, stating that non-zero net contributions from fluid stresses $\sigma n$ or additional external forces $F_e^i$ and torques $T_e^i$ lead to velocity variations.

The rotational speed $\angvel^i_k$ is linked with the orientation $R$ of the swimmer using

\begin{equation*}
	\begin{aligned}
		\frac{d}{dt} \Angle_k &= \angvel^i_k, \quad \text{for $k \in \{x,y,z\}$},\\
		\Rmat &=\Rmat_z(\Angle_z)\Rmat_y(\Angle_y)\Rmat_x(\Angle_x),
	\end{aligned}
\end{equation*}
where $R_k(\Angle_k)$ denotes the rotation matrix around axis $k \in \{x,y,z\}$ of 
angle $\Angle_k$. If $d=2$, $R(\theta)$ has the form: 
\begin{equation*}
	\begin{aligned}
		\Rmat(\Angle) &= \begin{bmatrix}
			\cos(\Angle) & \sin(\Angle)\\
			-\sin(\Angle) & \cos(\Angle)
		\end{bmatrix}.
	\end{aligned}
\end{equation*}
While in three dimensions:
$$	
	\Rmat_z(\Angle_z) = \begin{bmatrix}
		\cos(\Angle_z) & -\sin(\Angle_z) & 0\\
		\sin(\Angle_z) & \cos(\Angle_z)& 0 \\
		0& 0 & 1
	\end{bmatrix},
	\Rmat_y(\Angle_y) = \begin{bmatrix}
		\cos(\Angle_y) &0 & \sin(\Angle_y) \\
		0& 1& 0\\
		-\sin(\Angle_y) & 0& \cos(\Angle_y)
	\end{bmatrix},\\
$$
$$
	\Rmat_x(\Angle_x) = \begin{bmatrix}
		1 & 0 & 0\\
		0 &\cos(\Angle_x) & -\sin(\Angle_x) \\
		0 &\sin(\Angle_x) & \cos(\Angle_x) 
	\end{bmatrix}.
$$
The angle $\Angle$ belongs to  $\Theta$, with $\Theta = [-\pi,\pi]$ if $d=2$, 
or $\Angle_k \in \Theta$, $\Theta=[-\pi,\pi]\times[0,\pi]\times[0,\pi/2]$ if $d=3$.




\subsubsection{Collision model} 
\label{rigidcollision}

In this subsection, we describe how a swimmer's collisions with solid obstacles, other swimmers ,as active particles, and the boundaries of the fluid domain are handled. Our model is based on a short-range contact-avoidance repulsive force introduced by R. Glowinski in \cite{glowinski_fictitious_1998}. The force is activated in a collision zone of width $w_{col}$.
In equation \eqref{Eq:RB}, the external force and torque terms $F_e^i$ and $T_e^i$ summarize the physical interactions occurring when the distances $d_{ij}$ and $d_{i}$ between swimmer $i$ and another body $j$, or between swimmer $i$ and the domain boundary $\partial \Fluid_t$, are smaller than the width $w_{col}$ of the collision zone.  
The definition of the collision force for a swimmer-swimmer or swimmer-boundary pair is given by:


$$
F_{ij} = - \epsilon A  ( X_j - X_i)  1_{d_{ij} \leq w_{col}}, \qquad
F_{i \partial \Fluid_t} = - \epsilon' A ( X_{\partial \Fluid_t} - X_i )  1_{d_{i} \leq w_{col}}.
$$

Both equations contain an activation term $A$, depending on the type of the swimmer under consideration.  The vector connecting the contact points $X_i, X_j$ or $X_i, X_{\partial \Fluid_t}$
gives the direction of the force, and the stiffness parameters $\epsilon$ and $\epsilon'$ determine the force intensity. As a general rule, the force magnitude increases as the distance decreases. Finding the optimal values for $\epsilon$ and $\epsilon'$ is not trivial, since their values depend on fluid and swimmer properties.   

The total repulsion force applied on $\Solid^i$ is defined by:


$$
F_e^i = \sum_{(i,j)|d_{ij} \leq w_{col}} F_{ij} +  \sum_{i | d_i \leq w_{col}} F_{i \partial \Fluid_t},
$$

where the first sum runs over all the body pairs $(i,j)$ such that $d_{ij}\le w_{col}$, and the second sum runs over all body-boundary pairs such that $d_{i}\le w_{col}$.

Force $F_e^i$ leads to body rotations via its associated torque, defined by: 
$$
T_e^i = - (X_i - x_{CM}^i )\times F_e^i,
$$
where  $X_i$ is the contact point where $F_e^i$ is applied.

The total repulsion force and the external associated torque are added to the Newton equation \eqref{Eq:RB}, thus modifying the trajectory of the swimmer.  


\section{Numerical discretization}
\label{Sec:NumDiscr}

\subsection{The Arbitrary-Lagrangian-Eulerian formalism}
The fluid problem defined on a time-dependent domain is solved using the Arbitrary-Lagrangian-Eulerian (\ALE) formalism, which allows to follow the evolution of the fluid-structure interface.
In the rest of this subsection, we set the number of swimmers to $N=1$ to avoid the complexity of the notation.


Let $\CompDomain_t$ be the current computational domain where the fluid equations are solved. 
We define the \ALE\ maps $\Alemap^t:\CompDomain_0 \to \CompDomain_t$ as the family of smooth and bijective functions that describe the evolution of the computational domain. 
These functions are defined through the extension of the displacement field at the boundary of the swimmer $\partial \mathcal{S}_t$ to the interior of $\CompDomain_t$. 

For instance, if $\bar{\phi}^t:\partial \mathcal{S}_0 \to \partial \mathcal{S}_t$ is the boundary displacement, a possible definition of the \ALE\ maps is $\Alemap^t(X) = X + \phi^t(X)$, for $X \in \CompDomain_0$, where $\phi^t(X)$ is the extension of $\bar{\phi}^t$ via
\begin{equation}
	\left\{
	\begin{aligned}
		&\nabla \cdot((1+\tau(X))\nabla \phi^t(X)) = 0, \quad &\text{on $\CompDomain_0$},\\
		&\mathcal{\phi}^t= \bar{\phi}^t, &\text{on $\partial \mathcal{S}_0$},
	\end{aligned}
	\right.
 \label{Eq: VarFormALE}
\end{equation}
where $\tau(X)$ acts as a space-dependent diffusion coefficient influencing the regions where larger displacements are localised.

The time derivative of $\Vel$ in the \ALE\ frame can be expressed as a function of its Eulerian time derivative $\partial_t \Vel$ and the \ALE\ velocity $\Vel_{\mathcal{A}}(t,X)=\frac{\partial x}{\partial t}(t,(\Alemap^{t})^{-1}(x))$, $X \in \CompDomain_0$, which is the velocity of the moving domain:
\begin{equation}
\frac{\partial \Vel}{\partial t} \Big|_{\Alemap^t}(t,x)=\Vel_{\mathcal{A}} \cdot \nabla \Vel + \partial_t \Vel, \quad \text{$ x \in \CompDomain_t$}.
\end{equation}

In the discrete setting, it is not guaranteed that displacing the mesh according to the solution of \eqref{Eq: VarFormALE} always produces a valid triangulation: large displacements could lead to element inversions.
In order to prevent these problems, mesh quality measures \cite{field_qualitative_2000} are used to assess the validity of the triangulation: if the mesh deformation is ``small'' and the mesh quality remains above a predefined threshold, the domain is deformed according to the \ALE\ map; if the minimum of the mesh quality field falls below the threshold, the resulting mesh deformation is not viable and a remeshing procedure is applied before the \ALE\ map.

In our case, the discrete \ALE\ maps are computed by solving equation \eqref{Eq: VarFormALE} with piece-wise linear continuous finite elements. The spaces where the numerical solution and the test functions are chosen are defined as
\begin{equation}
	\begin{aligned}
		X_{\bar{\phi},h} &= \{  \phi , \phi \in [H^1(\CompDomain_{0})]^d\cap [\mathbb{P}_1(\CompDomain_{0})]^d, \phi=\bar{\phi} \text{ on $\partial \mathcal{S}_{0}$}\},\\
		X_{0,h} &= \{  \phi ,\phi \in [H^1_0(\CompDomain_{0})]^d\cap [\mathbb{P}_1(\CompDomain_{0})]^d\},
	\end{aligned}
\end{equation}
and the solution of the variational problem 
	\begin{equation}
		\begin{aligned}
			&\int_{\CompDomain_{0}} (1+\tau(X))\nabla \phi^{t_{n+1}}_h (X)  : \nabla v \,\dx= 0, \qquad &\text{$\forall v \in X_{0,h}$},
			\\
			&\mathcal{\phi}_h^{t_{n+1}}= \bar{\phi}^{t_{n+1}}, &\text{on $\partial\mathcal{S}_{0}$},
		\end{aligned}
  \label{eq:ALE_EDP}
	\end{equation}
where $A:B = \sum_{i,j} A_{ij}B_{ij}$ the Frobenius inner product, defines the new computational domain as $\CompDomain_{t_{n+1}}=\Alemap_h^{t_{n+1}}(\CompDomain_{0})$, where $\Alemap_h^{t_{n+1}}(X) = X + \phi_h^{t_{n+1}}(X)$. 

In the latter equations \eqref{eq:ALE_EDP}, $\tau$ is a piecewise constant coefficient, defined on each element $e$ of the domain's discretization as $ \tau\big |_e=(1-V_{min}/V_{max})/(V_e/V_{max})$, where $V_{max}$, $V_{min}$  and $V_e$ are the volumes of the largest, smallest and current element of the domain discretization \cite{kanchi_3d_2007}. This coefficient $\tau$ allows that the mesh deformation is  applied to elements of larger volume.

The evolution of the boundary of the swimmer is defined by both the rigid motion and its swimming gait by 
 
$$
\bar{\phi}^{t_{n+1}}(X) = \int_{0}^{t_{n+1}}\tvel + \angvel\times(X+\bar{\phi}^{t_n}(X)-X_{CM}- \phi^{t_n}(X_{CM})) + u_d(t,X)\, \mathrm{d}t
$$
where $X_{CM}$ is the mass center of the swimmer in the reference domain at initial time $t=0$. In others words, the center of mass at the current time in the current domain is defined by $x_{CM} := X_{CM} + \phi(X_{CM})$.

The time integration of $\bar{\phi}^{t_{n+1}}$  is performed in two steps: first, the contributions coming from the linear velocity are integrated to compute the new center of mass $\phi_1^{t_{n+1}}(X_{CM}) $ as
\begin{equation}
	\begin{aligned}
		\theta_{n+1}&= (t_{n+1}-t_n)\angvel^n + \theta_{n}, \quad R(\theta_{n+1}),\\
		\bar{\phi}_1^{t_{n+1}}(X) &= (t_{n+1}-t_n)\tvel^n + \bar{\phi}^{t_{n}}(X) + \int_{t_n}^{t_{n+1}} u_d(t,X)\, \mathrm{d}t;
	\end{aligned}
	\label{Eq:discDisp1}
\end{equation}
and then the orientation of the body is computed using the new center of mass 
\begin{multline}
	\bar{\phi}^{t_{n+1}}(X) \approx  R(\theta_{n+1}) (X+\bar{\phi}_1^{t_{n+1}}(X)- \phi_1^{t_{n+1}}(X_{CM})) + \phi_1^{t_{n+1}}(X_{CM})-X,
	\label{Eq:discDisp2}
\end{multline}

First, the rotation around the origin in the reference frame is performed, then the body is translated to its position in the current frame. A few fixed point iterations are performed at each time step to ensure the convergence of the body's position.

\subsection{Discretization of the fluid problem}

Let $\Vel_h^n$ and $\Pres_h^n$ denote the discrete approximations of the velocity and pressure fields at time $t_n$. Since the domain is time dependent, the functional spaces are time dependent as well via the discrete \ALE\ maps. In what follows, the functions $v$ and $p$ are defined in the current domain via the ALE map, whereas $\hat{v}$ and $\hat{p}$ are defined in the reference domain.
\begin{equation}
	\begin{aligned}
		V_h^t &= \{v: \CompDomain_t \to \mathbb{R}^d, \, v = \hat{v} \circ (\mathcal{A}_h^t)^{-1}, \, \hat{v} \in [H^1(\CompDomain_0)]^d \cap [\mathbb{P}_N(\CompDomain_0)]^d  \}, \\
		Q_h^t &= \{\Pres : \CompDomain_t \to \mathbb{R}, \, \Pres = \hat{\Pres} \circ (\Alemap_h^t)^{-1}, \, \hat{\Pres} \in\mathbb{P}_{N-1}(\CompDomain_0) \}.
	\end{aligned}
	\label{Eq:FEMspace}
\end{equation}

We choose the Taylor-Hood finite element spaces $\mathbb{P}_{2}$ for $V_h^t$ and $\mathbb{P}_{1}$ for $Q_h^t$. The discrete variational formulation of the Navier-Stokes equations in moving domain at time $t_{n+1}$ requires finding $(\Vel_h^{n+1}, \Pres_h^{n+1}) \in V_h^{t_{n+1}} \times Q_h^{t_{n+1}}$, $((\tvel^i)^{n+1},(\angvel^i)^{n+1})\in \mathbb{R}^d\times \mathbb{R}^{d^*}$ such that : 

\color{black}
\begin{eqnarray*}
\int_{\CompDomain_{t_{n+1}}} \rho_\Fluid \partial_t \Vel_h^{n+1}|_\Alemap  \cdot \tilde{\Vel} + \int_{\CompDomain_{t_{n+1}}} \rho_\Fluid ((\Vel_h^{n+1}-{\Vel_{\mathcal{A}}}_h^{n+1} )\cdot \nabla u_h^{n+1}) \cdot \tilde{\Vel}
 \\+ 2\mu \int_{\CompDomain_{t_{n+1}}} D(\Vel_h^{n+1}) : D( \tilde{\Vel}) -  \int_{\CompDomain_{t_{n+1}}} \Pres_h^{n+1} \nabla \cdot \tilde{\Vel} \\
 = \sum_{i=1}^N \int_{\partial \Solid_{t_{n+1}}^i} \sigma(\Vel_h^{n+1},\Pres_h^{n+1}) n \cdot \tilde{\Vel}
 + \int_{\CompDomain_{t_{n+1}}} f^{n+1} \cdot \tilde{\Vel},\\
 	\int_{\CompDomain_{t_{n+1}}} \tilde{\Pres} \nabla \cdot \Vel_h^{n+1} = 0,
  \end{eqnarray*}
  \begin{eqnarray*}
  \sum_{i=1}^N m^i \frac{d (\tvel^i)}{dt}^{n+1}\cdot \tilde{\tvel}^i = \sum_{i=1}^N \bigl[ F_e^i \cdot \tilde{\tvel}^i - \int_{\partial \Solid^i_{t_{n+1}}} \sigma \normal \cdot \tilde{\tvel}^i \bigr] ,
  \\ \sum_{i=1}^N \frac{d[R I^i R^T (\angvel^i)]^{n+1}}{dt} \cdot \tilde{\angvel}^i = \sum_{i=1}^N  \bigl[ T_e^i \cdot \tilde{\angvel}^i  -\int_{\partial \Solid^i_{t_{n+1}}} \sigma \times (x-\CenterMassi) \cdot  \tilde{\angvel}^i \bigr],
\label{Eq:discrete-var_form}
\end{eqnarray*}
for all $(\tilde{\Vel},\tilde{\Pres}) \in V_{h}^{t_{n+1}} \times Q_h^{t_{n+1}}, (\tilde{\tvel}^i,\tilde{\angvel}^i) \in \mathbb{R}^d\times \mathbb{R}^{d^*}$ for $i=1,\dots,N$.

We considered homogeneous Dirichlet and Neumann boundary conditions to simplify the formulations. In addition,
in the rest of this part, to state the discretization associated with these equations, we set the number of swimmers to $N=1$ to avoid the complexity of the notation.
Let us denote : 

\begin{eqnarray*}
a(u,v) := \int_{\CompDomain_{t}} \rho_\Fluid(\partial_t u|_\Alemap ) \cdot v + \int_{\CompDomain_{t}}\rho_\Fluid ((u -{u_{\mathcal{A}}}_h^{n+1} )\cdot \nabla u) \cdot v + 2\mu\int_{\CompDomain_{t}} D(u)	: D(v) \,, 
\end{eqnarray*}

$$
b(v,p) :=  -\int_{\CompDomain_{t}} p \nabla\cdot v \qquad \text{and} \qquad G(v) := \int_{{\CompDomain_{t}} } f \cdot v . 
$$
\color{black}

Following \cite{maury_direct_1999}, let us now denote the degrees of freedom that belong to the boundary of the swimmer by the subscript $\Gamma$ as $\Vel_\Gamma$ and the others by the subscript $I$ as $\Vel_I$.

\color{black}
We now discretize the previous equations as 

\begin{equation}
	\begin{bmatrix}
		A_{II} & A_{I\Gamma} & 0 & 0 & B_{I}^T  \\
		A_{\Gamma I} & A_{\Gamma \Gamma}  & 0 & 0 & B_{\Gamma}^T\\
		0 & 0 & T
		& 0 & 0 \\
		0 & 0 & 0& M
		& 0\\
		B_I & B_\Gamma & 0 & 0 & 0   \\
	\end{bmatrix}
	\begin{bmatrix}
		\Vel_I \\ \Vel_\Gamma \\ \tvel \\ \angvel \\ \Pres 
	\end{bmatrix} 
	= 
	\begin{bmatrix}
		G_I \\
		G_{\Gamma}\\
		F
		\\
		T
		\\
		0
	\end{bmatrix},
	\label{Eq:flu-rig-block}
\end{equation} 
\color{black}
where
\begin{align*}
&A_{JK} =	(a(\phi_{J_{i}},\phi_{K_{j}}))_{i,j},   \textrm{ for } J,K \in \{I,\Gamma\}  \textrm{ with $(\phi_{J_{i}})_i$ and $( \phi_{K_{j}} )_j$ is the basis of $V_h^{t_{n+1}}$},\\
&B_J = (b(\phi_{J_{i}},\psi_{j}))_{i,j} \textrm{ for } J \in \{I,\Gamma\} \textrm{ with $(\psi_j)_j$ is the basis of $Q_h^{t_{n+1}}$},\\
&G_J = G(\phi_{J_{i}}) \textrm{ for } J \in \{I,\Gamma\}, \quad T = m\mathbb{I}, \quad M = [RIR^T]^{n+1},\\
&F = F_e - \int_{\partial \Solid} \sigma \normal, \quad \textrm{and} \quad T = T_e -\int_{\partial \Solid} \sigma \times (x-\CenterMassi).
\end{align*}

In order to satisfy the boundary conditions $\Vel = \tvel+ \angvel \times (x-x_{CM}) + u_d$ on  $\partial \Solid$, we introduce the operator $\mathcal{P}$ such as 
\begin{equation}
\label{eq:boundary_constraint}
 (\Vel_I,\Vel_\Gamma,\tvel,\angvel,p)^T = \mathcal{P}	(\Vel_I,\tvel,\angvel,p)^T + (0,u_d\cdot\mathcal{A}_t^{-1},0,0,0)^T 
\end{equation}

with 
\begin{equation*}
	\mathcal{P} =
	\begin{bmatrix}
		\mathbb{I} &  0 & 0 & 0\\
		0 &\tilde{P}_\tvel & \tilde{P}_\angvel & 0\\
		0 &\mathbb{I} & 0 & 0\\
		0 &0 & \mathbb{I} & 0\\
		0 &0& 0& \mathbb{I}
	\end{bmatrix}\,.
\end{equation*}

In the previous matrix, $\tilde{P}_\tvel$ and $ \tilde{P}_\angvel$ are the interpolation operators that enable the expression of $u_\Gamma$ as a function of $\tvel$ and $\angvel$. 
Finally, plugging \eqref{eq:boundary_constraint} into \eqref{Eq:flu-rig-block}, we then obtain the discretization of the variational formulation of the fluid-swimmer problem 
\begin{equation*}
	\mathcal{P}^T \begin{bmatrix}
		A_{II} & A_{I\Gamma} & 0 & 0 & B_{I}^T  \\
		A_{\Gamma I} & A_{\Gamma \Gamma}  & 0 & 0 & B_{\Gamma}^T\\
		0 & 0 & T & 0 & 0 \\
		0 & 0 & 0& M 
		& 0\\
		B_I & B_\Gamma & 0 & 0 & 0   \\
	\end{bmatrix} \mathcal{P} \begin{bmatrix}
		\Vel_I \\ \tvel \\ \angvel \\ \Pres 
	\end{bmatrix}  
 = 
 \mathcal{P}^T 	
 \begin{bmatrix}
		G_I \\
		G_{\Gamma}\\
		F
		\\
		T
		\\
		0
	\end{bmatrix} - \mathcal{P}^T \begin{bmatrix}
		A_{II} & A_{I\Gamma} & 0 & 0 & B_{I}^T  \\
		A_{\Gamma I} & A_{\Gamma \Gamma}  & 0 & 0 & B_{\Gamma}^T\\
		0 & 0 & T & 0 & 0 \\
		0 & 0 & 0& M 
		& 0\\
		B_I & B_\Gamma & 0 & 0 & 0   \\
	\end{bmatrix}\begin{bmatrix}
		0 \\
		u_d \cdot \Alemap^{-1}_t\\
		0
		\\
		0
		\\
		0
	\end{bmatrix}
\end{equation*}.
\color{black}

\subsection{Collision model}

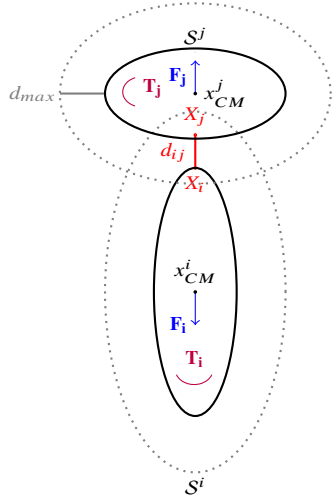
\begin{figure}
    
\begin{center}
\begin{tikzpicture}

\draw[color=red, thick](0,1.5*1.1) -- (0,1.9*1.1) node[font = \scriptsize, anchor = north east] {$d_{ij}$}  ;
\draw[color=red, thick](0,1.5*1.1) circle (0.25*1.1pt) node[font = \scriptsize, anchor = north] {$X_i$};
\draw[color=red, thick](0,1.9*1.1) circle (0.25*1.pt) node[font = \scriptsize, anchor = south] {$X_j$};

\draw[color=gray, thick](-1*1.2,2.4*1.1) -- (-1.5*1.2,2.4*1.1) ;

\node[font = \scriptsize,color=gray]  at (-1.8*1.2,2.4*1.1) {$d_{max}$};

\draw[color=blue, thin,->](0,0) -- (0,-0.4*1.1)  node[font = \scriptsize, anchor = east] {$\vec{F_i}$} ;
\draw[color=blue, thin,->](0,2.4*1.1) -- (0,2.8*1.1)node[font = \scriptsize, anchor = north east] {$\vec{F_j}$}  ;

\draw[color=black, thin](0,2.9*1.1) circle (0.0001) node[font = \scriptsize, anchor = south ] {$\Solid^j$} ;
\draw[color=black, thin](0,-2.8) circle (0.0001) node[font = \scriptsize, anchor = south] {$\Solid^i$} ;
\draw[color=purple, thin] (-0.8,2.85) arc(110:250:0.25 and 0.2)node[font = \scriptsize, anchor = south west] {$\vec{T_j}$} ;

\draw[color=purple, thin] (-0.25,-1.1) arc(200:340:0.25 and 0.2)node[font = \scriptsize, anchor = south east] {$\vec{T_i}$} ;

\draw[color=black, thick] (0,0) ellipse (0.5*1.1 and 1.5*1.1) ;
\draw[color=gray, thick, dotted] (0,0) ellipse (1*1.2 and 2*1.2);
\draw[color=black, thick](0,0.) circle (0.25*1.1pt) node[font = \scriptsize, anchor = south] {$x_{CM}^i$};

\draw[color=black, thick] (0,2.4*1.1) ellipse (1*1.2 and 0.5*1.2);
\draw[color=gray, thick, dotted] (0,2.4*1.1) ellipse (1.5*1.2 and 1.*1.2);
\draw[color=black, thick](0,2.4*1.1) circle (0.25*1.1pt) node[font = \scriptsize, anchor = west] {$x_{CM}^j$};

\end{tikzpicture}
\end{center}
   \caption{Notations for the collision model.  $\Solid^i,\Solid^j$ denote the solid objects, $d_{ij}$ the distance between the two bodies, and
$d_{max}$ the threshold distance for the narrow-band fast marching method. $F_i,F_j$ and $T_i, T_j$ denote the collision forces and torques.  }
    \label{fig:collision-tikz}

\end{figure}

The collision detection algorithm, which identifies the pairs of swimmers that are actually interacting as those whose surfaces are less than $w_{col}$ units apart, is based on the computation of distance functions. Figure \ref{fig:collision-tikz} shows the different notations. 
At each time instant, the distances $d_{ij}$ between the surfaces of two swimmers $\Solid^i$ and $\Solid^j$, $i \ne j$, are needed to determine if collision forces need to be activated. To compute these distances, a narrow-band variant of the fast marching method is used. 

When applied to a body $\Solid^i$, the  fast marching algorithm \cite{sethian_fast_1996} yields the distance field $D_i$ from $\Solid^i$ the to rest of the domain. However, since we are interested in the evaluation of the distance function in a small neighbourhood of $\Solid^i$, we employ a narrow-band approach, that only computes $D_i$ close to $\partial \Solid^i$. This choice accelerates considerably the computations, especially in three dimensions, and is also suitable for parallel execution.

The size of the neighborhood is set to a predefined threshold $d_{max}$, defined as a function of the width of the collision zone, $d_{max} \approx 2 w_{col}$. Upon reaching the threshold distance $d_{max}$, the narrow-band 
approach assigns a default value $\delta = d_{max}$ to the distance field, corresponding to the maximum value reached:
\begin{center}
	\begin{equation*}
		D_i^{NB} = \left\{\begin{array}{rcl}
		0 \;, \quad &\mbox{on }& \partial \Solid^i , \\
		D_i \;, \quad &\mbox{for }& D_i \leq d_{max} ,\\
		\delta\;, \quad &\mbox{elsewhere}&   .
		\end{array}\right.\;
	\end{equation*}
\end{center}

Using the two distance fields $D_i^{NB}$ and $D_j^{NB}$, the distance $d_{ij}$ is described by:
$$
d_{ij} = ||\arg\min_{x \in \partial \Solid^j}  D_i^{NB}(x) - \arg\min_{x \in \partial \Solid^i} D_j^{NB}(x)||_2. 
$$
The boundary points $X_i = \arg\min_{x \in \partial \Solid^i} D_j^{NB}(x)$ and $X_j = \arg\min_{x \in \partial \Solid^j}  D_i^{NB}(x)$ 
give the contact points of the swimmers $\Solid^i$ and $\Solid^j$, i.e., the coordinates 
of the points where they will interact. 

Similarly to the case of two interacting bodies, it is possible to compute the distance between a swimmer and the boundaries of the fluid domain: 
$$
d_{i} = ||\arg\min_{x \in \partial \Fluid_t} D_i^{NB}(x) - \arg\min_{x \in \partial \Solid^i} D_{\Fluid_t}^{NB}(x)||_2 \mbox{,}
$$
where $D_{\Fluid_t}^{NB}$ represents the distance field obtained by applying the fast marching method 
to the portion of the domain boundary $\partial \Fluid_t \backslash ( \cup_i \partial \Solid^i ) $.

\section{Applications}
\label{Sec:Applications}
%

We present numerical simulations of swimming micro-organisms to illustrate the framework we have presented in this chapter. The results are obtained using the \feelpp\ finite element library \cite{prudhomme_feelppfeelpp_2023}.

\subsection{Flagellated swimmer: 3D sperm cell}
The previous framework allows the simulation of flagellated swimmers, provided that the analytical expression of the deformation velocity $\Vel_d(x,t)$ is known in advance. This is the case for a sperm cell which propagates planar waves along its flagellum, where $\Vel_d(x,t)$ has the following form \cite{razavi_ale-based_2015}
\begin{equation}
	\Vel_d(t,X)=
	\begin{bmatrix}
		\frac{2\pi}{4T}[A_{max}(X-X_j)/L]^2\frac{2\pi}{\lambda}\cos(4\pi (\frac{t}{T}-\frac{X}{\lambda})) \\
		\frac{2\pi}{T}[A_{max}(X-X_j)/L]\cos(2\pi (\frac{t}{T}-\frac{X}{\lambda}))
	\end{bmatrix}.
	\label{Eq:up}
\end{equation}
Equation \eqref{Eq:up} represents the velocity of a sinusoidal wave of linearly increasing amplitude, that propagates from the head of the swimmer to the tip of its flagellum. 
While the $y$-component of $\Vel_d$ comes from the time derivative of the sinusoidal wave
\begin{equation*}
	Y(t,X) = \frac{A_{max}(X-X_j)}{L} \sin\Bigg(2\pi \Bigg(\frac{t}{T}-\frac{X}{\lambda}\Bigg)\Bigg),
\end{equation*}
the $x$-component ensures the non-extensibility of the sperm tail \cite{taylor_analysis_1951}. 
In the previous equations, $T$ is the period of the wave, $\lambda$ is its the wavelength, $L$ is the length of the flagellum, $X_j$ is the coordinate of the head-flagellum junction and $A_{max}$ is the amplitude at the flagellum's distal end.

\begin{figure}
	\centering
	\includegraphics[width=0.45\linewidth]{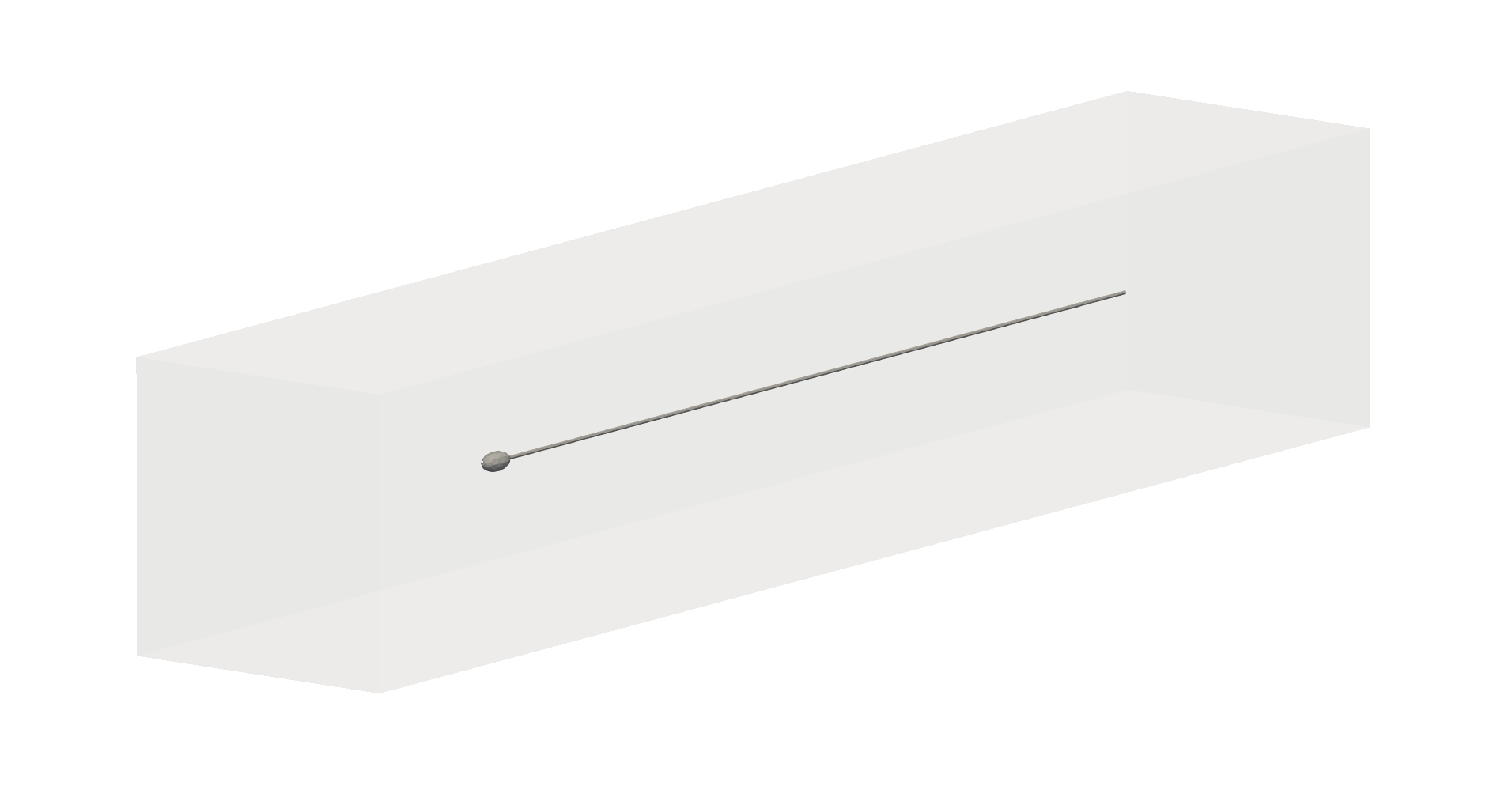}
	\includegraphics[width=0.45\linewidth]{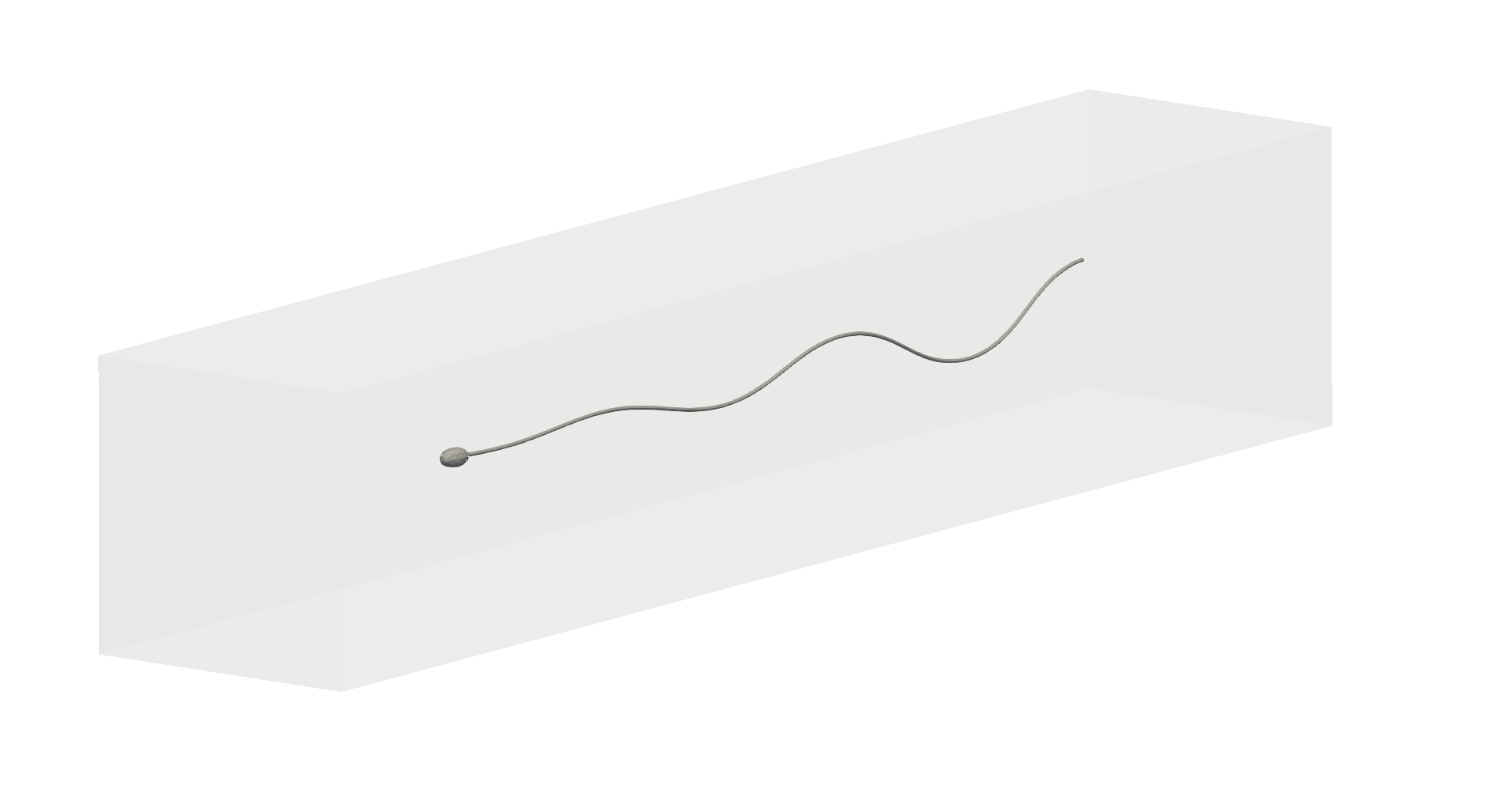}
	\includegraphics[width=0.45\linewidth]{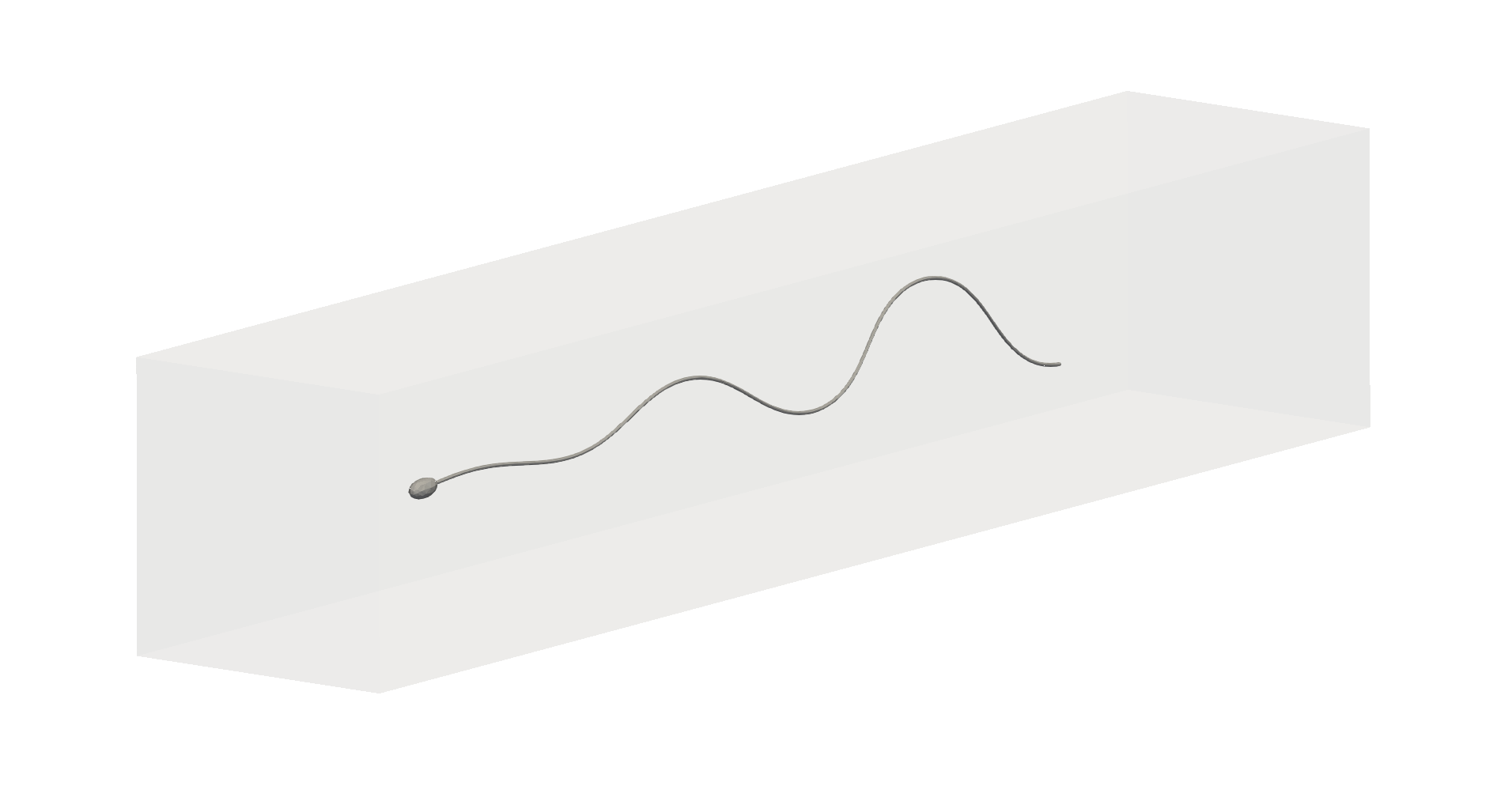}
	\includegraphics[width=0.45\linewidth]{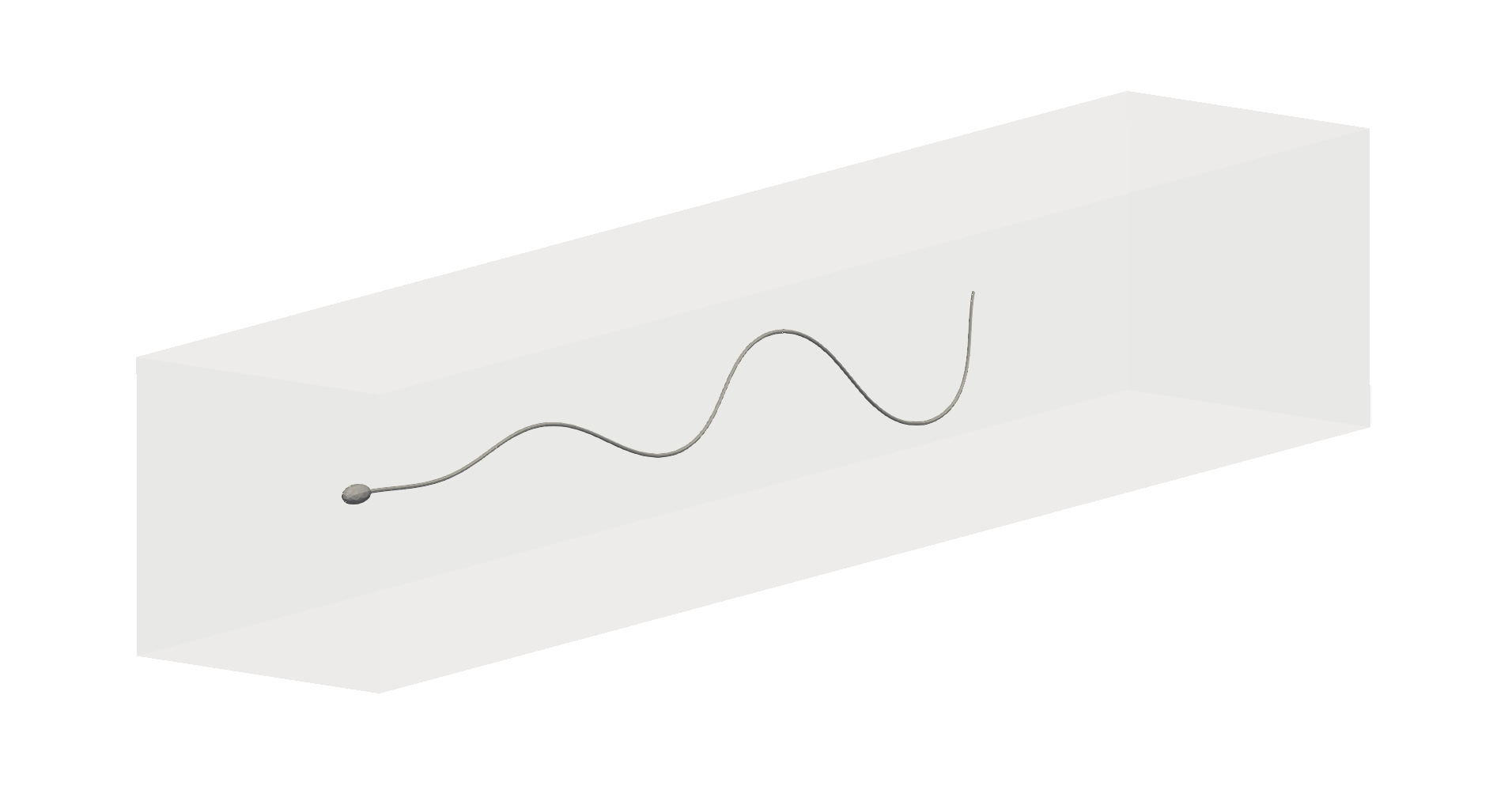}
	\caption{Three dimensional simulation of a swimming sperm cell propagating planar waves on its flagellum. Flagellar beating is restricted to $xy$ plane, and the maximal amplitude of the wave is $A_{max}=\SI{4}{\micro\meter}$. The four images show the position and shape of the sperm cell at the time instants $t=\SI{0}{\second}$, $t=\SI{0.5}{\second}$, $t=\SI{1.85}{\second}$ and $t=\SI{2.65}{\second}$. }
	\label{Fig:3dsperm_planar}
\end{figure}

In Figure \ref{Fig:3dsperm_planar} we show the position and shape of the sperm cell at four time instants, which were determined by solving \eqref{Eq:NS}-\eqref{Eq:RB} by prescribing $u_d$ as in \eqref{Eq:up} and by parameterizing the equation as in \cite{razavi_ale-based_2015}. In particular, the wave is restricted to the $xy$ plane and its maximal amplitude is $A_{max} = \SI{4}{\micro\meter}$. Moreover, the deformation velocity $\Vel_d$ is reached by gradually increasing $A_{max}$ in time until reaching $A_{max}=4$. We obtain that the swimmer moves on a straight line with a constant speed, once the wave is fully developed.

To validate the fluid-structure interaction model and ALE framework, that we introduced in this chapter, we compared the observed swimming velocity of a $2$D spermatozoon to the results presented in \cite{razavi_ale-based_2015}. Good agreements are found. 

\subsection{Multi-body swimmer: the three-sphere swimmer}

\begin{figure}
	\centering
	\includegraphics[scale=0.55]{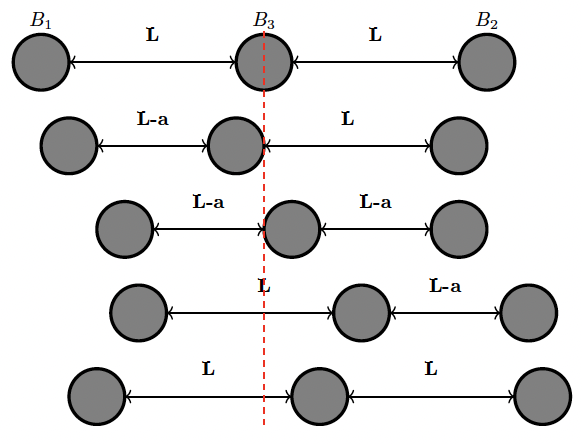}

	\caption{Representation of the three-sphere swimmer and its swimming gait. The gait is composed of four strokes in which one of the rods is alternatively retracted or elongated.}
	\label{Fig:3ss1}
\end{figure}

In this section, we consider the well-known three-sphere swimmer \cite{najafi_simple_2004}, a model swimmer which is extensively utilized in micro-swimming due to its simplicity and its capability to capture complex hydrodynamic effects.
This type of swimmer consists of three identically-sized spheres connected by rods which are alternatively extended and retracted to produce a net motion. The sequence begins with retracting the left rod, followed by the right one. Then, the left rod is extended to reach its initial length and finally, the right rod too (see Figure \ref{Fig:3ss1}). 
This sequence of four movements results in a straight motion when the swimmer is far from boundaries. 

In this case, we solved the Stokes equations by imposing $\rho_{\Fluid}=0$ in the system \eqref{Eq:NS}-\eqref{Eq:RB} and by defining the deformation velocity $u_d$ as a function of the relative speeds between the spheres. The ALE map is not affected by this change form Navier-Stokes to Stokes equations. More details are given in \cite{berti_modelling_2021}. 

\begin{figure}[h!]
    \centering
    \centerline{\includegraphics[scale=0.5]{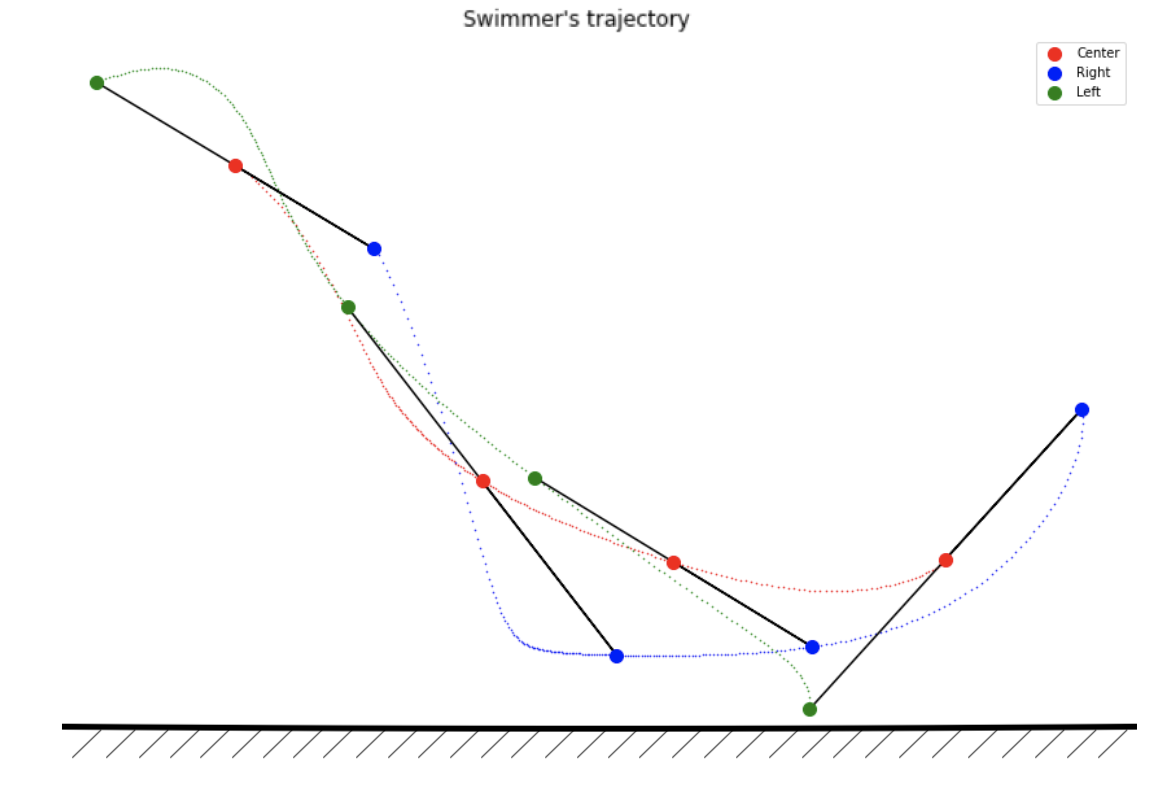}}
    \caption{Behavior of a three-sphere swimmer swimming towards the boundary of 
	the computational domain. The swimmer reorients due to contact forces and hydrodynamic interactions with the plane wall.}
    \label{fig:3ss}
\end{figure}

\begin{figure}[h!]
    \centering
    \includegraphics[scale=0.175]{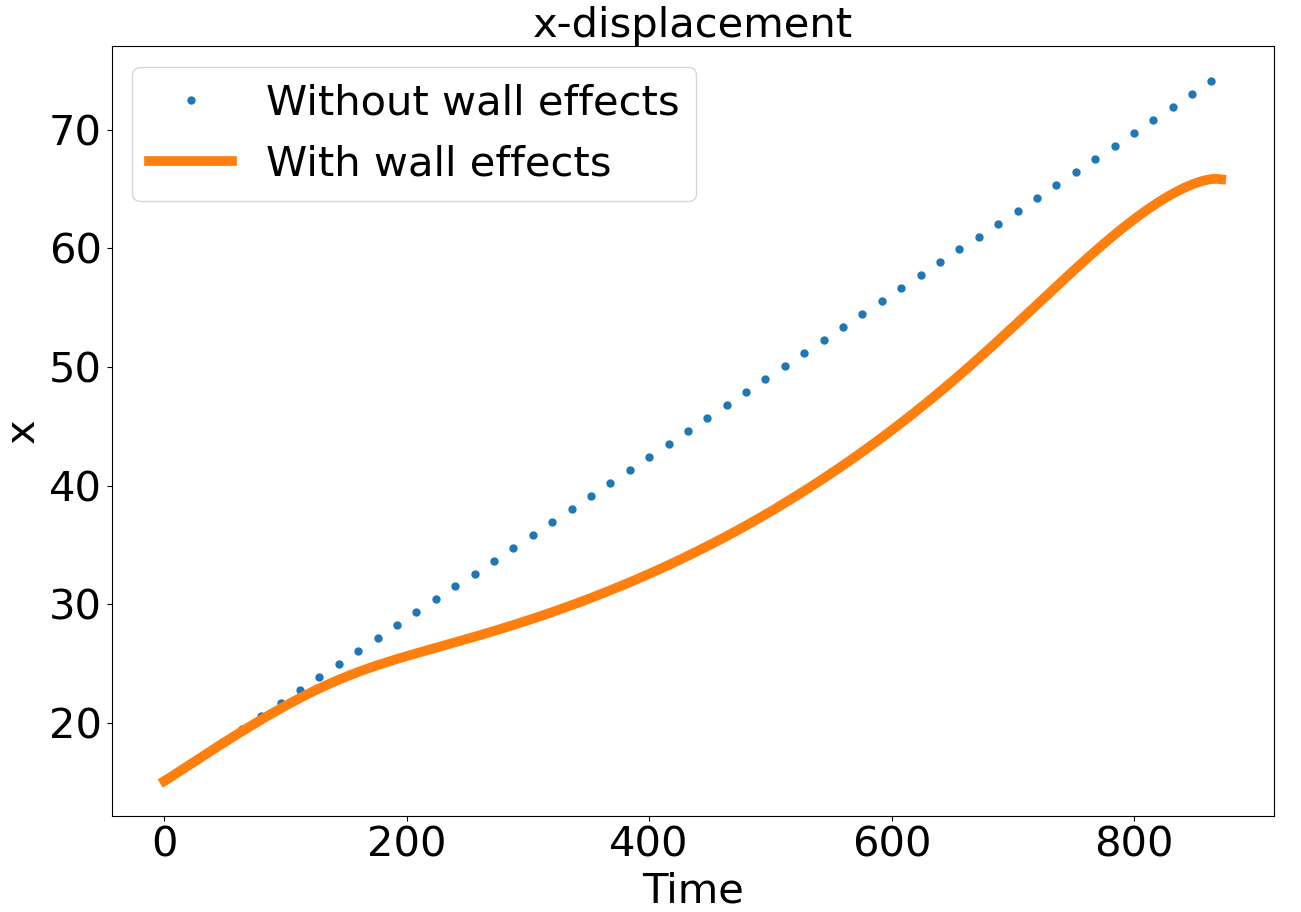}
    \includegraphics[scale=0.175]{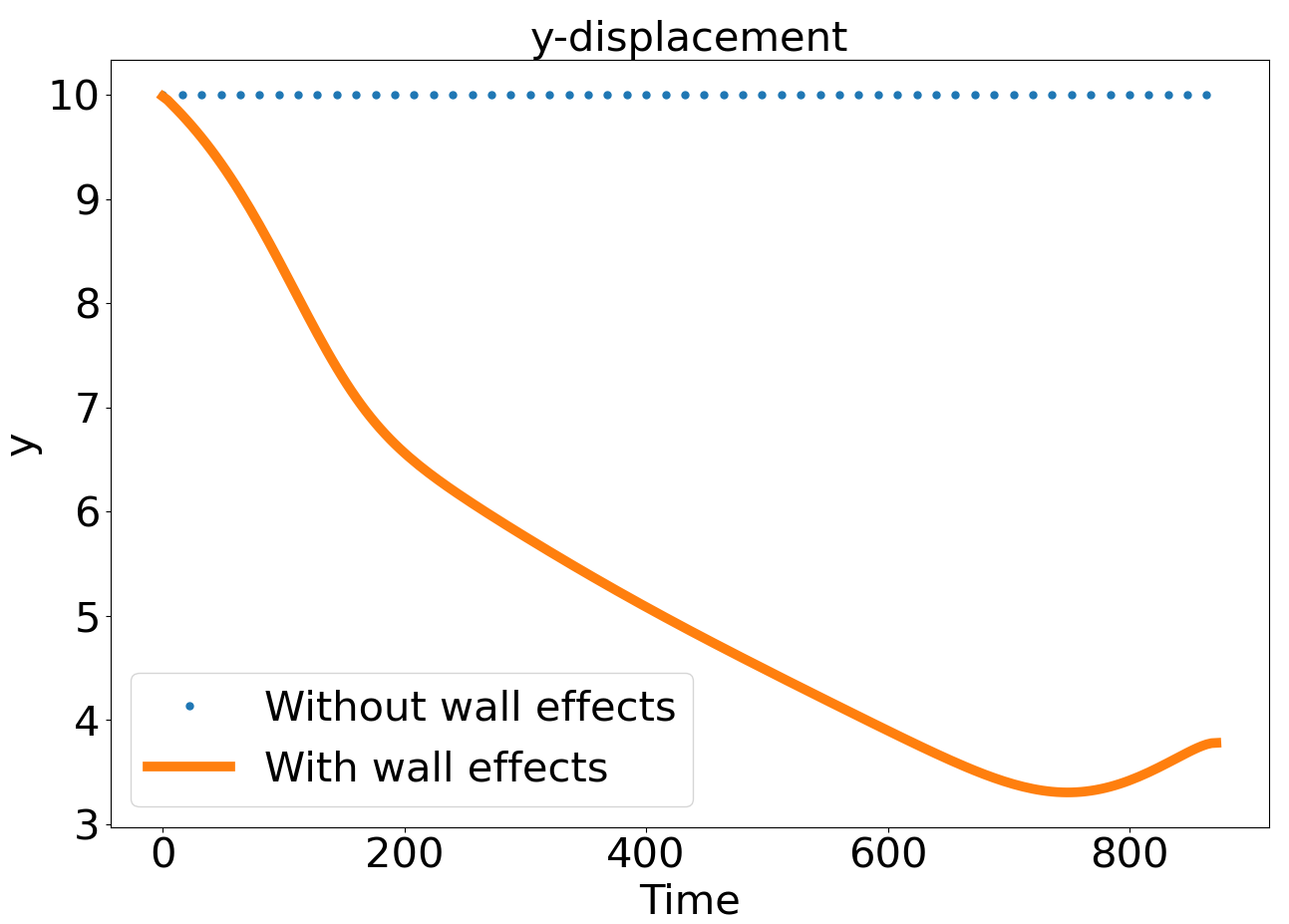}
    \caption{Examples of vertical and horizontal displacements of three-sphere swimmers. The orange continuous lines corresponds to the $x$ and $y$ trajectories of a swimmer heading towards the boundary of the domain and reorienting after the contact with the wall. The blue dotted line describe the trajectories of a swimmer swimming far from the wall and parallel to it.}
    \label{fig:3ss_2}
\end{figure}

Figures \ref{fig:3ss} and \ref{fig:3ss_2} show how the motion of the swimmer is affected by the presence of a plane wall. Figure \ref{fig:3ss} captures the behavior of the three-sphere swimmer as it approaches and interacts with the boundary, while in Figure \ref{fig:3ss_2} two behaviours are presented. First, the orange continuous lines describes the trajectory of a swimmer that, due to its initial orientation and swimming strategy, gets closer to the boundary of the channel. Once its right sphere arrives in the collision zone, collision forces are applied, and the swimmer starts to change direction. It continues rotating until its left sphere reaches the collision zone, where the repulsive force acting on this sphere propels the swimmer upwards, distancing it from the boundary. 
Secondly, the blue dotted lines correspond to the trajectory of a swimmer whose rods are parallel to the plane wall and which is not perturbed by its presence. The displacement of this latter swimmer is in good agreement with the literature \cite{najafi_simple_2004}.

\subsection{Rigid bodies with tangential velocities: squirmers}

\begin{figure}
	\centering

        \includegraphics[scale=0.2]{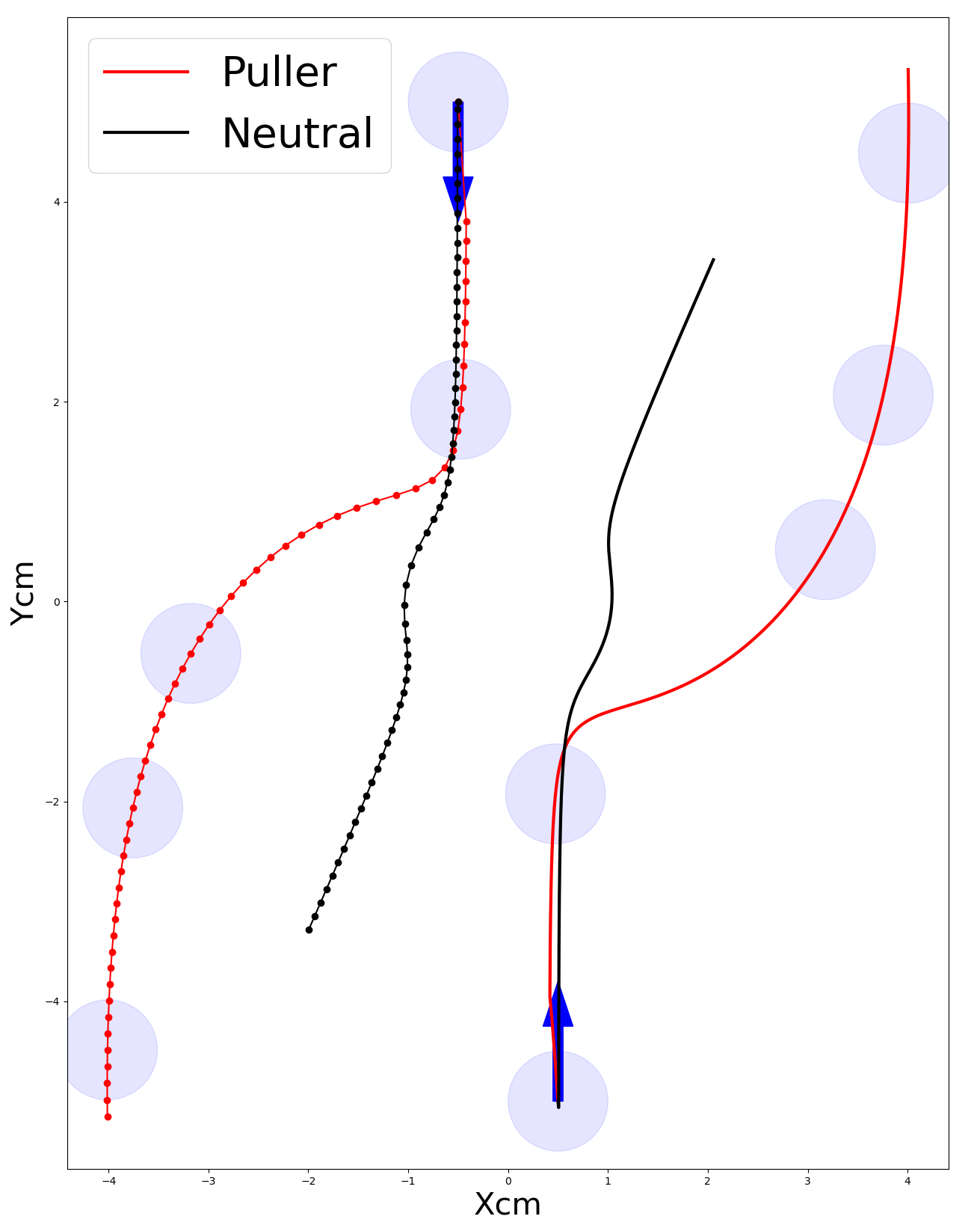}
    
	\includegraphics[scale=0.14]{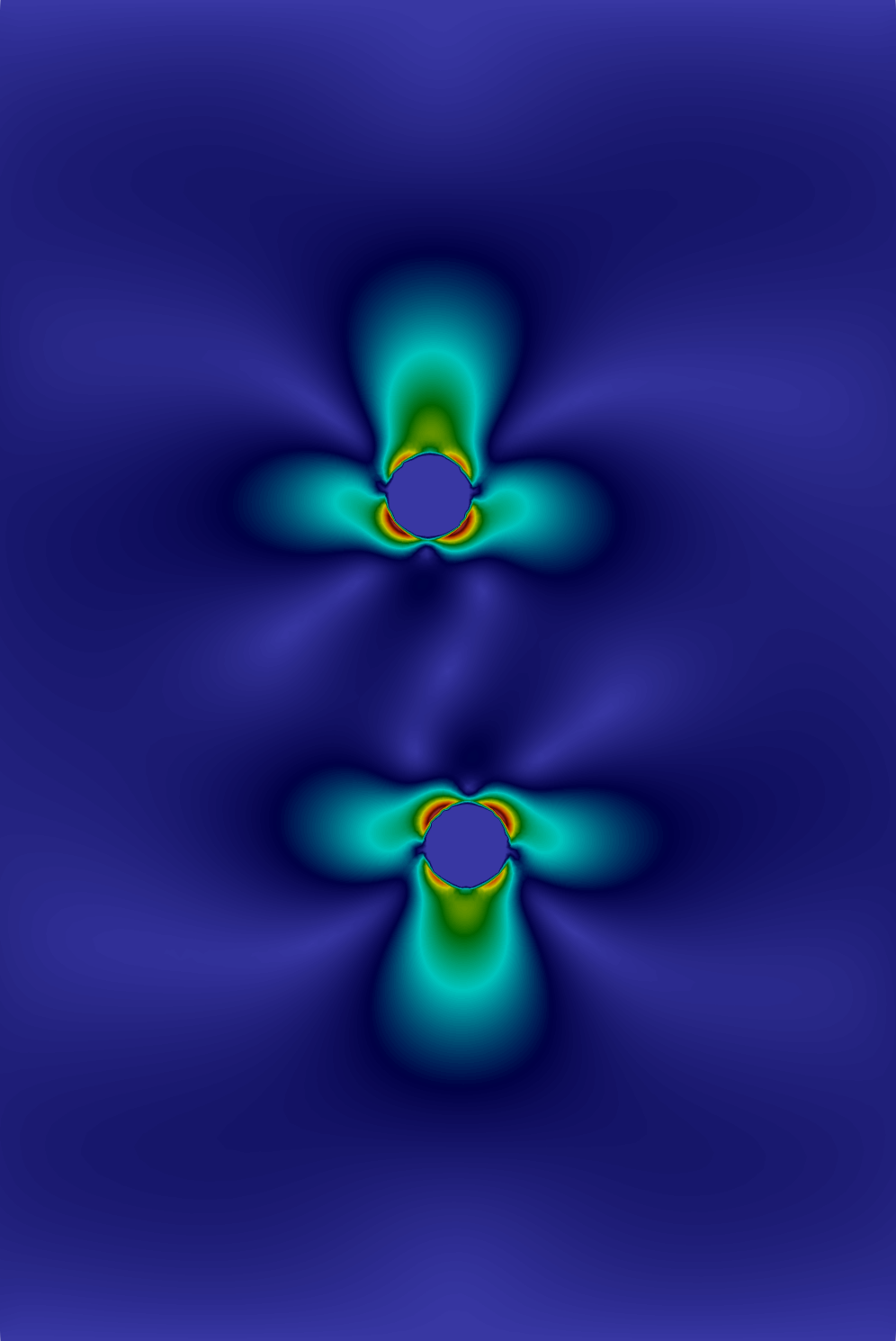}
	\includegraphics[scale=0.14]{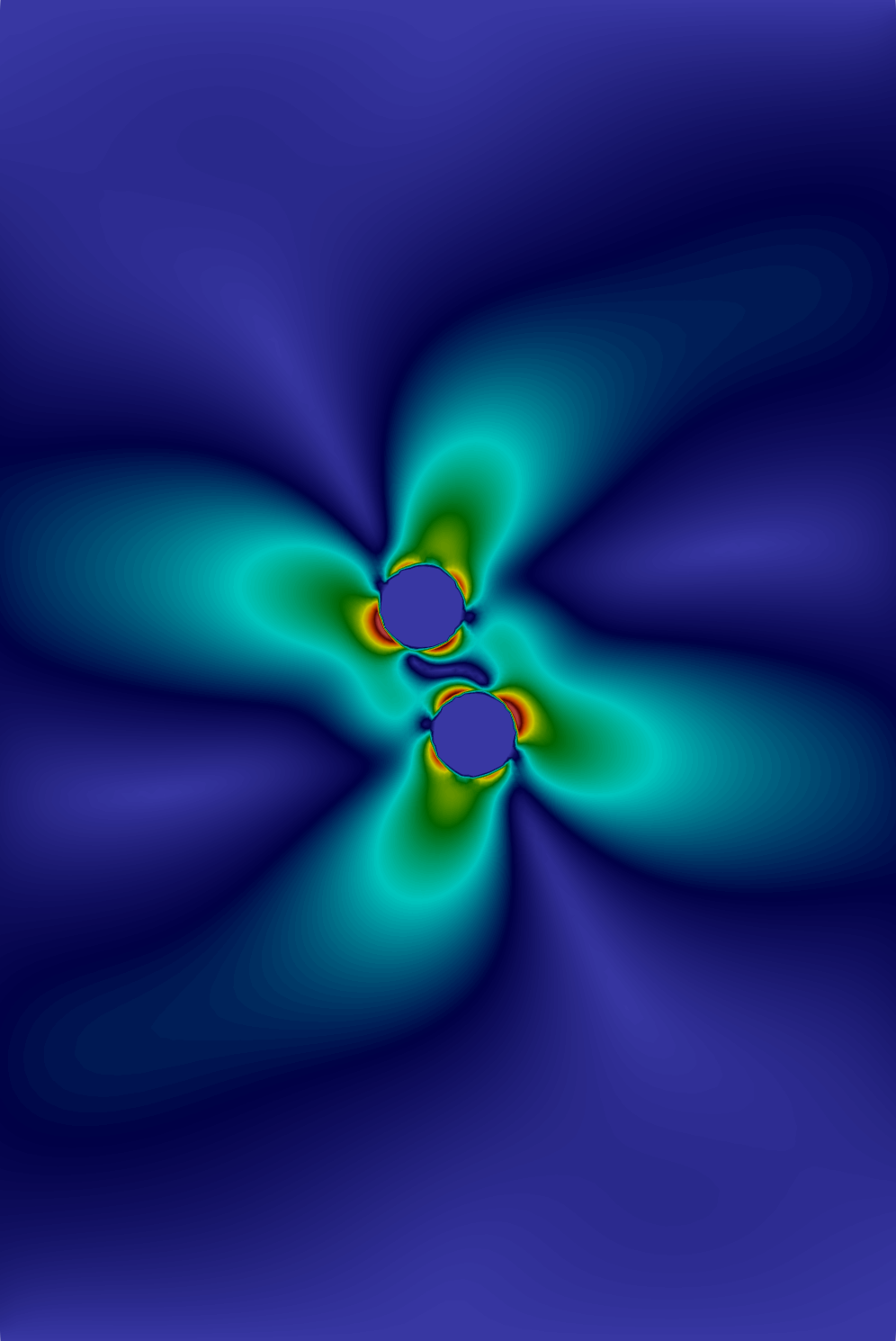}
	\includegraphics[scale=0.14]{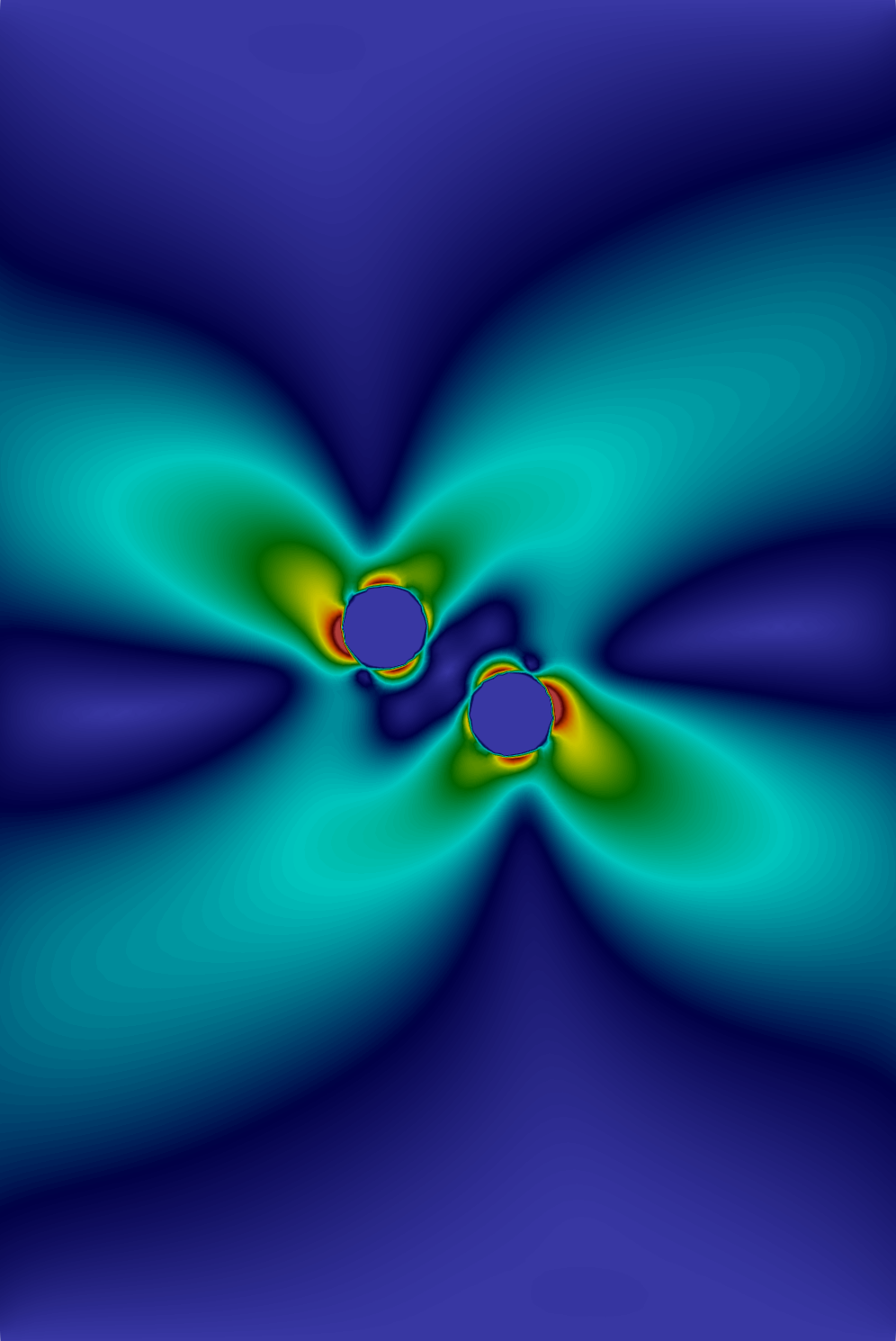}
	\includegraphics[scale=0.14]{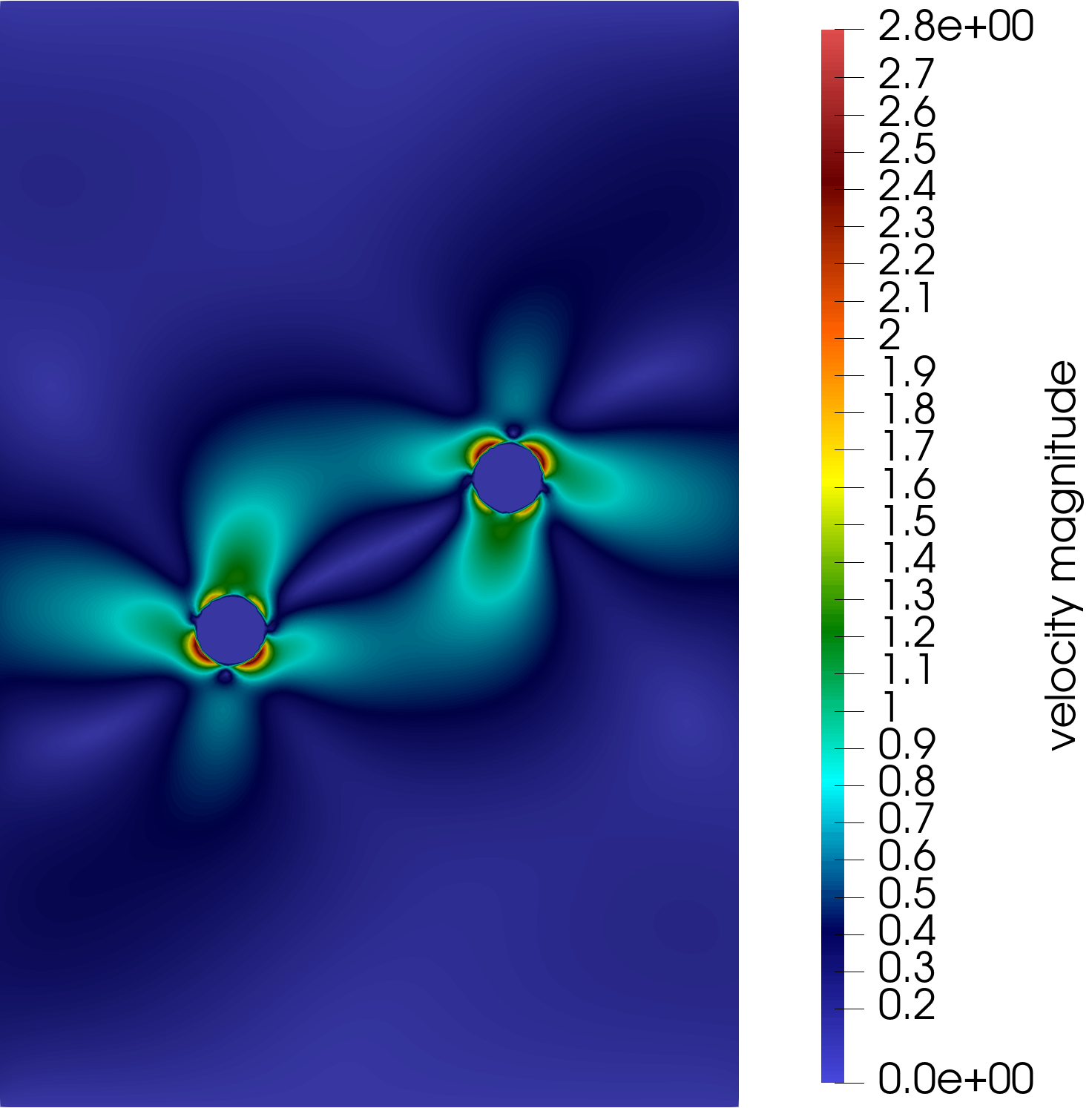}

	\caption{In the top figure: trajectories of two interacting squirmers. Starting from the same initial configuration, the positions of two neutral squirmers (black lines) and pullers (red lines) are shown as they swim towards one another. The dotted lines show the trajectory of the upper squirmers and the smooth lines the one of the lower squirmers. The two arrows correspond to their initial orientation. Depending on the propulsion type, the trajectories show different behaviours.
 In the bottom figures: positions of the squirmers and magnitude of the fluid velocity are shown at different time instants of the interaction.}
	\label{Fig:squirmers}
\end{figure}

Some ciliated micro-organisms can be approximated via the squirmer model \cite{blake_spherical_1971,lighthill_squirming_1952}, which considers them to be rigid bodies with prescribed velocity patterns at the surface. Most frequently, the swimming gait is encoded in the function $u_d$ by prescribing a time-independent velocity field tangent to $\partial \Solid_0^i$, modeling the cilia waving pattern on the surface of the body. For a circular swimmer moving in direction $\overrightarrow{e}$, this velocity field is prescribed by:

$$
u_d(x) =  B_1 \Bigg[ 1 + \beta (\overrightarrow{e} \cdot \overrightarrow{r}) \Bigg] \Bigg[(\overrightarrow{e} \cdot \overrightarrow{r})\overrightarrow{r} - \overrightarrow{e} \Bigg]
$$ 
where  $ \overrightarrow{r} = \frac{x - x_{CM}} {||x - x_{CM}||} $, $B_1$ and $\beta$ are the swimming speed and propulsion type ( $\beta < 0$ corresponds to a pusher, $\beta > 0$ to a puller and $\beta = 0$ to a neutral squirmer).\\
In order to showcase the usage of the collision algorithm between swimmers, we simulate the interaction between two squirmers of the same propulsion type, considering neutral squirmers or pullers. The trajectories in Figure \ref{Fig:squirmers} show that, the two swimmers change orientation due to the collision forces when getting closer, and then move away from each other.
The intensity of the repulsion depends on the type of squirmers: in the case we have considered, neutral squirmers reach a smaller distance than pullers before deviating from their initial trajectory. Similar behaviours for squirmer-squirmer interactions are found in \cite{ishikawa_hydrodynamic_2006}.

\subsection{Motion of a collection of solids: inside the zebrafish arteries}

\begin{figure}
	\centering
	\includegraphics[width=0.85\linewidth]{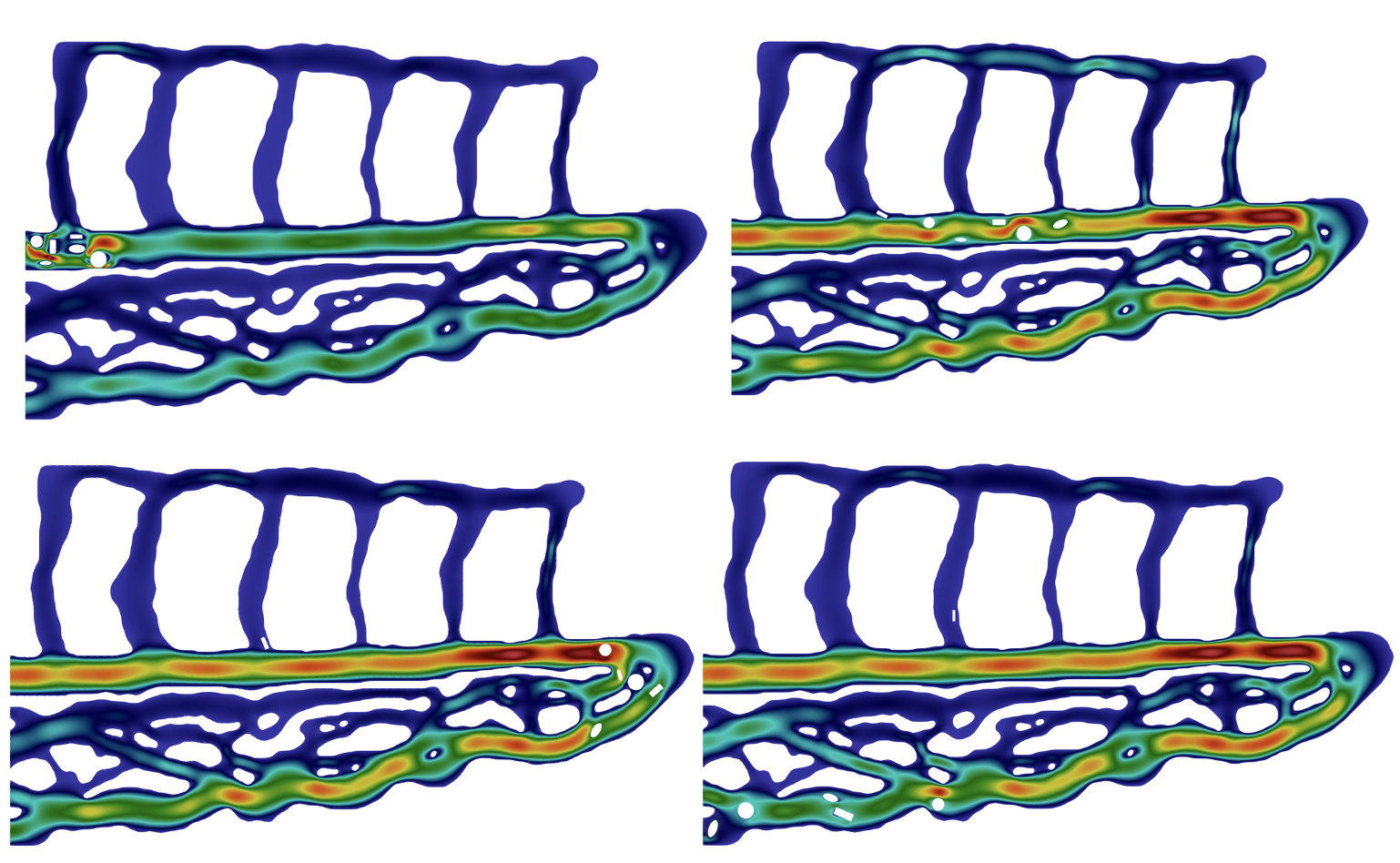}

	\caption{Moving solid bodies in a complex two-dimensional geometry. Due to collision forces and torques, the bodies start to rotate and take different trajectories to cross the zebrafish.}
	\label{Fig:insertion2D}
\end{figure}

\color{black}
In biological processes, such as particle transport in blood vessels, solid bodies move in geometrically complex domains. The purpose of this last application is to demonstrate that our numerical framework allows simulating such a biological phenomenon. The trajectory of a collection of solids, within a two-dimensional reconstruction of the vascular system of a zebrafish, a model used in cancer biology, is shown. 
\color{black}
The solids are initially placed to the inlet boundary of the arterial network and their motion is driven by a pulsatile velocity imposed at this boundary, given by:  
$$
u = 35 * |\sin \bigl( \frac{\pi}{0.15}*t \bigr) |.
$$ 
In this application, remeshing 
is necessary, since mesh deformation via ALE maps alone is not sufficient to guarantee the good quality of the mesh.
The snapshots in Figure \ref{Fig:insertion2D} show the positions at different time instants of the solids moving in the complex geometry of the zebrafish. Depending on their shape and initial position, the objects have different trajectories within the network. Additionally, interactions with boundaries and other objects lead to rotational motion, causing the solids to be pushed from the main artery into regions with lower fluid velocity.




\section{Conclusion}

\color{black}
In this paper, we present a framework using finite element methods with the Arbitrary Lagrangian-Eulerian (ALE) approach to simulate the dynamics of rigid deformable swimmers immersed into a Navier-Stokes fluid. Our simulations account for collision effects in both swimmer-swimmer interactions and interactions between swimmers and walls. Additionally, our approach allows the computation of swimmer dynamics in complex geometrical environments. All implementations are carried out using the finite elements library Feel++ \cite{prudhomme_feelppfeelpp_2023}. Several numerical examples are provided, showing good agreement with the literature.
\color{black}



\bibliographystyle{plain}
\bibliography{references}

\end{document}